\documentclass[11pt]{article}
\usepackage{fullpage}
\usepackage{latexsym}
\usepackage{epic}
\usepackage{eepic}
\usepackage{graphics}

\newcommand{\cmt}[1]{\ifhmode\newline\fi{\sf *** \ \ #1 \\}}

\newtheorem{theorem}{Theorem}
\newtheorem{lemma}[theorem]{Lemma}
\newtheorem{corollary}[theorem]{Corollary}

\newtheorem{prop}[theorem]{Proposition}

\newcommand\proofend{\qed}
\newenvironment{proof}{\noindent\textbf{Proof:}}{\qed\hfill\\\mbox{}\par}

\newcommand{\qed}{\mbox{\rule{7pt}{7pt}}}

\newcommand{\eps}{\varepsilon}

\newcommand{\TT}{{\cal T}}

\newcommand{\R}{{\mathbf{R}}}

\newcommand{\N}{{\mathbf{N}}}

\newcommand{\heading}[1]{\vspace{1ex}\par\noindent{\bf #1}}

\def\:{\colon}

\newcommand\ExpOp{\mathbf{E}}
\newcommand{\Ex}[1]{\ExpOp\left[#1\right]}
\newcommand{\Var}[1]{\mathrm{Var}\!\left[#1\right]}

\newcommand{\Prob}[1]{\mathbf{P}\!\left[#1\right]}

\long\def\onefigure#1#2{
\begin{figure*}[tbh]
\begin{center}
#1
\end{center}
\caption{#2}
\end{figure*}
} 

\newcommand{\lipefig}[2]  
{\onefigure{\def\IPEfile{match#1.ipe}\input{\IPEfile}}{\label{f:#1} #2} }

\newcommand\asize{r}
\newcommand\bsize{s}
\newcommand\anode{a}
\newcommand\bnode{b}
\newcommand\aclass{A}
\newcommand\bclass{B}
\newcommand\nblks{q}
\newcommand{\lis}{\mathrm{LIS}}

\newcommand{\CMM}{Ctr.~de Modelamiento Matem\'{a}tico}

\newcommand{\DIM}{Dept.~Ing.~Matem\'{a}tica}

\newcommand{\UChile}{U.~Chile}

\begin{document}

\title{Expected length of the longest common subsequence for large alphabets}

\author{
	{\sc Marcos Kiwi}\thanks{Gratefully acknowledges the support of
	ICM P01--05 and Fondecyt 1010689.}
	\\
	  {\footnotesize Depto.~Ing.~Matem\'{a}tica and}\\[-1.5mm]
	  {\footnotesize Ctr.~Modelamiento Matem\'atico UMR 2071, }\\[-1.5mm]
	  {\footnotesize University of Chile}\\[-1.5mm]
	  {\footnotesize Correo 3, Santiago 170--3, Chile}
	\\[-1.5mm]{\footnotesize e-mail: \texttt{mkiwi@dim.uchile.cl}}
	\and
	{\sc Martin Loebl}\thanks{Gratefully acknowledges the support of ICM-P01-05.
This work was done while visiting the \DIM, \UChile.}
	\\
	   {\footnotesize Dept.~of Applied Mathematics and}\\[-1.5mm]
	   {\footnotesize Institute of Theoretical Computer Science (ITI)}\\[-1.5mm]
	   {\footnotesize  Charles University}\\[-1.5mm]
	   {\footnotesize  Malostransk\'{e} n\'{a}m. 25, 118~00~~Praha~1}\\[-1.5mm]
	   {\footnotesize Czech Republic}
	\\[-1.5mm]{\footnotesize e-mail: \texttt{loebl@kam.mff.cuni.cz}}
	\and
	{\sc Ji\v{r}\'{\i} Matou\v{s}ek}\thanks{This research was done
	while visiting the \CMM, UMR--UChile 2071, \UChile,
        Santiago, supported by Fondap in Applied
	  Mathematics 2000--05.}
	\\
	   {\footnotesize Dept.~of Applied Mathematics and}\\[-1.5mm]
	   {\footnotesize Institute of Theoretical Computer Science (ITI)}\\[-1.5mm]
	   {\footnotesize  Charles University}\\[-1.5mm]
	   {\footnotesize  Malostransk\'{e} n\'{a}m. 25, 118~00~~Praha~1}\\[-1.5mm]
	   {\footnotesize Czech Republic}
	\\[-1.5mm]{\footnotesize e-mail: \texttt{matousek@kam.mff.cuni.cz}}
	     }
	\date{}

	\maketitle

\begin{abstract}
We consider the length $L$ of the longest common subsequence
of two randomly uniformly and independently chosen $n$ character
words over a $k$-ary alphabet.
Subadditivity arguments yield that $\Ex{L}/n$ converges to a constant
$\gamma_k$.
We prove a conjecture of Sankoff and Mainville from the early 80's
claiming that $\gamma_k\sqrt{k}\to 2$ as $k\to\infty$.

\end{abstract}

\section{Introduction}\label{s:intro}
Consider two sequences of length $n$, with
  letters from a size $k$ alphabet $\Sigma$, say $\mu$ and $\nu$.
The longest common subsequence (LCS) problem is that of finding
the largest value $L$ for which there are
$1\leq i_1< i_2< \ldots < i_L \leq n$ and
 $1\leq j_1< j_2< \ldots < j_L \leq n$
 such that $\mu_{i_t} = \nu_{j_t}$, for all $t=1,2,\ldots,L$.

The LCS problem has emerged more or less independently in several remarkably
  disparate areas, including the comparison of versions of computer
  programs, cryptographic snooping,
  and molecular biology.
The biological motivation of the problem is that long molecules such
  as proteins and nucleic acids like DNA
  can be schematically represented as sequences from a finite alphabet.
Taking an evolutionary point of view, it is natural to compare two
  DNA sequences by finding their closest common ancestors.
If one assumes that these molecules evolve only through the process
  of inserting new symbols in the representing strings, then ancestors
  are substrings of the string that represent the molecule.
Thus, the length of the longest common subsequence of two strings is a
  reasonable measure of how close both strings are.
In the mid 1970's, Chv\'atal and Sankoff~\cite{cs75}
proved that the expected
  length of the LCS of
  two random $k$-ary sequences of length $n$ when normalized
  by $n$ converges to a constant.
The value of this constant $\gamma_k$
  is unknown although much effort has been spent in
  finding good upper an lower bounds for it (see, for example, \cite{bgns99}
and references therein).
The best known upper and lower bounds for $\gamma_k$ do not have a
 closed form.
There are  obtained either as numeric approximation to the
  solutions of a nonlinear equation or as a
  numeric evaluation of some series expansion (see~\cite{dancikthesis}
  for a survey of such results).

Although the problem of determining $\gamma_{k}$ has a simple statement,
  it has turned out to be a challenging mathematical endeavor.
Moreover, its quite naturally motivated.
Indeed, a claim that two DNA sequences of length $n$ are
  far apart makes sense provided their LCS differs significantly
  from $\gamma_{4}n$ (since DNA sequence have $4$ basis elements).

We analyze the behavior of $\gamma_{k}$ for $k$
tending to infinity, and more generally, we
consider the expected length of the LCS when $k$
is an (arbitrarily slowly growing) function of $n$ and $n\to\infty$.
The focus on the case where $k$ grows with $n$ is partly inspired by
 the work of Kiwi and Loebl~\cite{kl02}.
For a bipartite graph $G$ over two size $n$
  totally ordered color classes $A$ and $B$, they considered
\[
L(G) = \max\{L:
  \mbox{$\exists a_1<\ldots< a_L$, $b_1<\ldots< b_L$,
   $a_ib_i\in E(G)$, $1\leq i\leq L$}\}\,,
\]
and studied its behavior when $G$ is uniformly chosen among
 all possible $d$-regular bipartite graphs on $A$ and $B$.
They established that $L_{n}(G)/\sqrt{dn}\to 2$ as $n\to\infty$
  provided $d=o(n^{1/4})$.
Under this latter condition, any node of the $d$-regular
  bipartite graph can potentially be matched to a $d/n\to 0$
  fraction of the other color class nodes.
In the case of interest here, that is the LCS problem with
  $k\to\infty$,
  it also happens that any sequences' character can be matched
  to an expected $1/k\to 0$ fraction of the other sequence's characters.
Both for this work and in~\cite{kl02}, the vanishing fraction
  of (expected) potential matches is a key issue.

In this paper we confirm a conjecture of
  Sankoff and Mainville from the early 80's~\cite{sm83} stating that
\begin{equation}\label{eqn:lim}
\lim_{k\to\infty} \gamma_{k}\sqrt{k} = 2\,.
\end{equation}
(See~\cite[\S~6.8]{pevzner00} 
  for a discussion of work on lower and upper bounds on $\gamma_k$
  as well as a stronger version, due to Arratia and Steele, of the 
  above stated conjecture.)

The constant $2$ in~(\ref{eqn:lim}) arises from a connection
  with another celebrated problem known as the longest increasing
  sequence (LIS) problem.
The problem is also referred to as ``Ulam's problem.'' 
  (e.g., in~\cite{kingman73,bdj99,okounkov00}).
Some (e.g.,~\cite{pevzner00}) incorrectly credit Ulam for raising it in~\cite{ulam61}
  where he mentions (without reference) 
  a ``well--known theorem'' asserting that given
  $n^{2}+1$ integers in any order, it is always possible to find among 
  them a monotone subsequence of $n+1$.
The theorem is due to  Erd\H{o}s and Szekeres~\cite{es35}.
The discussion in~\cite{ulam61} concerns only the behavior
  of the monotonic subsequence of a randomly and uniformly 
  chosen permutation of $n^{2}+1$ elements.
Monte Carlo simulations are reported in~\cite{bb67}, where it is 
  observed that over the range $n\leq 100$, the limit of the LIS
  of $n^2+1$ randomly chosen elements, 
  when normalized by $n$, approaches $2$.
Hammersley~\cite{hammersley72} gave a rigorous proof of the existence of the limit
  and conjectured it was equal to $2$.
Later, Logan and Shepp~\cite{ls77}, based on a result by
  Schensted~\cite{schensted61}, proved that $\gamma\geq 2$;
  finally, Vershik and Kerov~\cite{vk77} obtained that $\gamma\leq 2$.
In a major recent breakthrough due to Baik, Deift, Johansson~\cite{bdj99}
  the asymptotic distribution of the longest increasing sequence
  random variable has been determined.
For a detailed account of these results, history and related work
  see the surveys of Aldous and Diaconis~\cite{ad99} and
  Stanley~\cite{stanley02}.

\medskip
It has been speculated that the behavior of the longest
  strictly/weakly increasing subsequence of a uniform
  random word of length $n$, with letters from $\Sigma$
  may have ``connections with the subject of sequence comparison
  statistics, motivated by DNA sequence matching \ldots''~\cite{ad99}.
Our work re-enforces this speculation and in fact does more.
It partly elicits the nature of the connection and the conditions under
  which sequence matching statistics relate to the behavior of
  longest increasing sequences.

\section{Statement of Results}\label{s:contrib}
Let $A$ and $B$ henceforth denote two disjoint totally
 ordered sets.
We assume that the elements of $A$ are numbered $1,2,\ldots,|A|$
  and those of $V$ are numbered $1,2,\ldots,|B|$.
We denote by $\asize$ and $\bsize$ the size of $|A|$ and $|B|$,
 respectively.
Typically, we have $\asize=\bsize=n$.

Now, let $G$ be a bipartite graph with color classes $A$ and $B$.
Two distinct edges $ab$ and $a'b'$ of $G$ are said to be \emph{noncrossing}
  if $a$ and $a'$ are in the same order as $b$ and $b'$; in other
  words, if $a<a'$ and $b<b'$ or $a'<a$ and $b'<b$.
A matching of $G$ is called \emph{planar} if every distinct pair
  of its edges is noncrossing.
We let $L(G)$ denote the number of edges of a maximum
 size planar matching in $G$ (note that $L(G)$ depends on
 the graph $G$ \emph{and} on the ordering of its color classes).

We will focus on the following two models of random graphs:
\begin{itemize}
\item[--] The {\em random words model\/}
$\Sigma(K_{n,n};k)$: the distribution over the set of
subgraphs of $K_{n,n}$  obtained by uniformly and independently
assigning each  node of $K_{n,n}$
one of $k$ characters and keeping those edges whose endpoints are
associated to equal characters. Note that only
disjoint unions of complete bipartite graphs may appear in this model.

\item[--] The {\em binomial random graph model\/}
$G(K_{n,n};p)$: the
distribution over the set of subgraphs of $K_{n,n}$ where each
edge of $K_{n,n}$ is included with probability $p$, and these
events are mutually independent. (This is an obvious modification
of the usual $G(n,p)$ model for bipartite graphs with ordered
color classes.)
\end{itemize}

In order to keep the presentation simple, we first formulate
and prove the results  for the random words model.
Then, in Section~\ref{s:bin}, we state analogous results for
the binomial random graph model.
These results' proofs are almost identical to the case of the
random words model, and we only briefly comment on them.

Our results essentially say that $L(\Sigma(K_{n,n};k))\cdot \sqrt k/n$
converges to $2$ as $k\to\infty$, provided that $n$ is sufficiently
large in terms of $k$.

\begin{theorem}\label{t:}
For every $\eps>0$ there exist $k_0$ and $C$ such that for all
$k>k_0$ and all $n$ with $n/\sqrt k>C$ we have
\[
     (1-\eps)\cdot \frac{2n}{\sqrt k} \ \leq\ 
	     \Ex{L(\Sigma(K_{n,n};k))}\ \leq\  (1+\eps)\cdot \frac{2n}
	{\sqrt k}\,.
\]
Moreover, there is an exponentially small tail bound; namely, for every
$\eps>0$ there exists $c>0$ such that for $k$ and $n$ as above,
$$
\Prob{\left| L(\Sigma(K_{n,n};k))- \frac{2n}
        {\sqrt k}\right|\geq \eps \frac{2n}
	        {\sqrt k}}\leq e^{-cn/\sqrt k}.
$$
\end{theorem}

\begin{corollary}
The limit
  $\gamma_k=\lim_{n\to\infty} \Ex{L(\Sigma(K_{n,n};k))/n}$ exists,
  and
\[
\lim_{k\to\infty} \gamma_k\sqrt{k} = 2.
\]
\end{corollary}

\section{Tools}\label{s:tools}
The crucial ingredient in our proofs 
is a sufficiently precise result on 
the distribution of the length of the longest increasing
  subsequence in a random permutation.
We state a remarkable strong result of 
 Baik, Deift and Johansson \cite[eqn.~(1.7) and~(1.8)]{bdj99}
 (our formulation slightly weaker than theirs,
 in order to make the statement simpler).
A much weaker tail bound than provided by them
  would actually suffice for our proof.

\begin{theorem}\label{t:bdj}
Let $\lis_{N}$ be the random variable corresponding to the
length of the longest increasing subsequence
of a randomly chosen permutation of
$\{1,\ldots,N\}$. There are positive constants $B_0,B_1$, and $c$
such that  for every $\lambda$ with $B_0/N^{1/3}\leq \lambda\leq \sqrt{N}-2$,
 \[
\Prob{\lis_N\geq 2\sqrt N + \lambda\sqrt{N}}
\leq B_1 \exp\left(-c\lambda^{3/5}N^{1/5}\right), 
\]
and for every $\lambda$ with $B_0/N^{1/3}\leq \lambda\leq 2$,
 \[
    \Prob{\lis_N \leq 2\sqrt{N} -\lambda \sqrt{N}}
      \leq B_1 \exp\left(-c\lambda^{3}N\right).
 \]
\end{theorem}

We will also need a suitable version of Talagrand's inequality;
see, e.g., 
 \cite[Theorem 2.29]{jlr00}.

\begin{theorem}[Talagrand's inequality]\label{t:tala}
Suppose that
$Z_1,\ldots,Z_N$ are independent random variables taking their values
  in some set $\Lambda$.
Let $X=f(Z_1,\ldots,Z_N)$, where $f:\Lambda^N\to \R$
  is a function such that the following two conditions hold for some
  number $c$ and a function $\psi$:
\begin{itemize}
\item[(L)] If $z,z'\in \Lambda^N$ differ only in the $k$th coordinate,
  then $|f(z)-f(z')|\leq c$.
\item[(W)] If $z\in \Lambda^N$ and $r\in\R$ with $f(z)\geq r$, then there
  exists a witness $(\omega_j:j\in J)$,
  $J\subseteq\{1,\ldots,N\}$, $|J|\leq \psi(r)/c^{2}$, such that for all
  $y\in\Lambda^N$ with $y_i=\omega_i$ when $i\in J$, we have $f(y)\geq r$.
\end{itemize}
Let $m$ be a median of $X$.
Then, for all $t\geq 0$,
\[
\Prob{X\geq m+t} \leq 2e^{-t^2/4\psi(m+t)}.
\]
and
\[
\Prob{X\leq m-t} \leq 2e^{-t^2/4\psi(m)}.
\]
\end{theorem}

We will also need the following version of Chebyshev's inequality:
\begin{lemma}\label{l:cheb}
Let $X_1,\ldots,X_N$ be  random variables attaining values $0$ and $1$,
and let $X=\sum_{i=1}^N X_i$.
Let $\Delta=\sum_{i\neq j} \Ex{X_iX_j}$.
Then, for all $t> 0$,
\[
\Prob{\,\left|X-\Ex{X}\right| \geq t}
  \leq {1\over t^{2}}\Bigl(\Ex{X}(1-\Ex{X})+\Delta\Bigr)\,.
\]
\end{lemma}
\begin{proof}
Since $\Prob{\left|X-\Ex{X}\right| \geq t} \leq \Var{X}/t^2$
  and
\begin{eqnarray*}
\Var{X}	  & = & \sum_{i,j}(\Ex{X_iX_j}-\Ex{X_i}\Ex{X_j}) \\
	  & = & \sum_{i}\Ex{X^2_i}-\sum_{i,j}\Ex{X_i}\Ex{X_j}
		+\sum_{i\neq j}\Ex{X_iX_j}\,,
\end{eqnarray*}
the desired conclusion follows by additivity of expectation and
  the fact that since $X_i$ is
  an indicator variable, $X^2_i=X_i$.
\end{proof}

\section{Small graphs}\label{s:small}
In this section we  derive a result essentially 
saying that Theorem~\ref{t:} holds if $k$ is sufficiently
  large in terms of $n$.
For technical reasons, we also need to consider bipartite
 graphs with color classes of unequal sizes.

\begin{prop}\label{p:small}
For every $\delta>0$, there exists a (large) positive constant $C$
such that:
\begin{enumerate}
\item[\rm (i)]
If $\asize\bsize\geq Ck$ and 
$(\asize+\bsize)\sqrt{\asize\bsize}\leq \delta k^{3/2}/6$,
then  with $m_u=m_u(\asize,\bsize)=2(1+\delta)\sqrt{\asize\bsize/k}$, we have
\[
 \Prob{L(\Sigma(K_{\asize,\bsize};k))\geq m_u+t}\leq 2e^{-t^2/8(m_u+t)}
\]
for all $t\geq 0$.
\item[\rm (ii)]
If $\asize\bsize\geq Ck$ and 
and $\asize+\bsize\leq \delta k/6$, then with
$m_u$ as above and $m_l=m_l(\asize,\bsize)=2(1-\delta)\sqrt{\asize\bsize/k}$,
we have 
\[
 \Prob{L(\Sigma(K_{\asize,\bsize};k))\leq m_l-t}\leq 2e^{-t^2/8m_u}
\]
for all $t\geq 0$.
\end{enumerate}
\end{prop}

Let $G$ be a random bipartite graph generated according
to the random words model $\Sigma(K_{\asize,\bsize};k)$.
The idea of the proof is simple: we show that (ignoring 
degree $0$ nodes) $G$ is ``almost''
a matching, and the size of the largest planar matching in a
random matching corresponds precisely to the length of the
longest increasing sequence in a random permutation of the
appropriate size. 

We have to deal with the (usually few) vertices
of degree larger than one. To this end, we define 
a graph $G'$ obtained from $G$ by removing all
edges incident to nodes of degree at least $2$.
Throughout, $E$ and $E'$ denote $E(G)$ and $E(G')$, respectively.

We clearly have $\Ex{|E|}=\asize\bsize/k$. We will need a tail bound
for large deviation from the expectation; a simple second-moment
argument (Chebyshev's inequality) suffices.

\begin{lemma}\label{l:num-edges-up}
For every $\eta>0$, 
\[
  \Prob{ \left| |E|- {\asize\bsize\over k}\right|
  \geq \eta \cdot {\asize\bsize\over k}} \leq
    {1\over \eta^2(\asize\bsize/k)}\,.
\]
\end{lemma}
\begin{proof}
For $e\in E(K_{\asize,\bsize})$ let $X_e$ be the indicator of the
  event $e\in E$.
Furthermore, let $X=|E|=\sum_{e\in E} X_e$.
The $X_{e}$'s are indicator random variables with expectation
  $1/k$. 
Moreover, since $\Ex{X_eX_f}=1/k^{2}$ for $e\neq f$, we have
 $\sum_{e\neq f}\Ex{X_eX_f}=\asize\bsize(\asize\bsize-1)/k^{2}
				 =(\Ex{X})^{2}-\Ex{X}/k$.
Thus, Lemma~\ref{l:cheb} yields
\[
\Prob{|X-\Ex{X}|\geq \eta\Ex{X}}
  \leq {1\over \eta^{2}\Ex{X}}\left(1 - {1\over k}\right).
\]
The desired conclusion follows immediately.
\end{proof}

Now we bound above the expectation of $|E\setminus E'|$.
\begin{lemma}\label{l:exp-del-edges}
\[
\Ex{|E\setminus E'|}\leq (\asize+\bsize){\asize\bsize\over k^2}\,.
\]
\end{lemma}
\begin{proof}
Let $Y_w$ equal the degree of $w$ if it is at least $2$ and $0$ otherwise.
Define $Y=\sum_{w\in V(G)} Y_w$.
Note that $|E\setminus E'|\leq Y$ (equality does not necessarily
  hold since both endpoints of an edge might be incident on nodes of degree
  at least $2$).
Let $P_d$ be the probability that
  a vertex in color class $\aclass$ has exactly $d$ incident edges.
For any node $\anode$ in color class $\aclass$,
\[
\Ex{Y_\anode} = \sum_{d=2}^{\bsize}dP_d=
\Ex{\deg_G(\anode)}-P_1
   = {\bsize\over k}-{\bsize\over k}\left(1-{1\over k}\right)^{\bsize-1}
  \leq \left({\bsize\over k}\right)^{2}
\]
(using $(1-x)^h\geq 1-hx$).
Similarly $\Ex{Y_\bnode}\leq (\asize/k)^{2}$ for all
  nodes $\bnode$ in color class $\bclass$,
  and so
\[
\Ex{|E\setminus E'|}\leq\Ex{Y}\leq (\asize+\bsize){\asize\bsize\over k^2}\,.
\]
\qed

\heading{Proof of Proposition~\ref{p:small}. }
Changing one of the characters associated to a vertex of
  a bipartite graph $G$ changes the value of $L(G)$ by at most $1$.
Hence $L(G)$ is $1$-Lipschitz.
Furthermore, the characters associated to $2\omega$ nodes of $G$
  suffice to certify the existence of $\omega$ noncrossing edges (and thus
  $L(G)\geq \omega$).
So Talagrand's inequality applies and, with $m$ denoting
  a median of $L(G)$, yields
\[
\Prob{L(G)\geq m+t}\leq 2e^{-t^2/8(m+t)}\ \ \mbox{and}\ \ 
\Prob{L(G)\leq m-t}\leq 2e^{-t^2/8m}.
\]
The proposition will follow once we show that $m_l\leq m\leq m_u$.
To prove that $m\leq m_u$, it suffices to verify that
\begin{equation}\label{ineq:median-bound}
 \Prob{L(G)\geq m_u} \leq {1\over 2}\,.
\end{equation}
Let $\eta>0$ be a suitable real
  parameter which we will specify later.
We observe that
  since $|E'|\leq |E|$ and $L(G)-L(G')\leq |E\setminus E'|$,
\begin{eqnarray*}
\Prob{L(G)\geq m_u} &\leq &
    \Prob{|E|\geq (1+\eta){\asize\bsize\over k}} \\
 & & \mbox{}+ \Prob{|E\setminus E'|\geq\delta\sqrt{\asize\bsize\over k}\,}\\
 & & \mbox{}+ \Prob{L(G')> (2+\delta)\sqrt{\asize\bsize\over k},
	    |E'| < (1+\eta){\asize\bsize\over k}}.
\end{eqnarray*}

We bound the terms one by one.
By Lemma~\ref{l:exp-del-edges} and Markov's inequality,
\begin{eqnarray}\label{ineq:del-edges-up}
\Prob{|E\setminus E'|\geq \delta\sqrt{\asize\bsize\over k}\,}
  \leq {\asize+\bsize\over \delta k}\sqrt{\asize\bsize\over k}
  \leq {1\over 6}\,.
\end{eqnarray}

Taking $N=(1+\eta)\asize\bsize/k$ and
  $\lambda=[(2+\delta)/\sqrt{1+\eta}\,]-2>0$
  in Theorem~\ref{t:bdj}, we get that
\begin{eqnarray}\label{ineq:bdj}
\lefteqn{\Prob{L(G')\geq (2+\delta)\sqrt{\asize\bsize\over k},
  |E'|< (1+\eta){\asize\bsize\over k}}} \\ \nonumber
   & & \mbox{} \leq B_1\exp\left(-c\lambda^{3/5}N^{1/5}\right)
   \leq B_1\exp\left(-c\lambda^{3/5}\left({rs\over k}\right)^{1/5}\right)\,.
\end{eqnarray}
\relax From Lemma~\ref{l:num-edges-up},
 (\ref{ineq:del-edges-up}),
  and (\ref{ineq:bdj}), it follows that
\begin{eqnarray*}
\Prob{L(G)\geq m_u}
  \leq {1\over \eta^2 (rs/k)} + {1\over 6}
   + B_1\exp\left(-c\lambda^{3/5}\left({\asize\bsize\over k}\right)^{1/5}\right)\,.
\end{eqnarray*}
So, (\ref{ineq:median-bound}) follows by taking, say, $\eta=\sqrt{6/C}$
and using $\asize\bsize \geq Ck$.

To establish that $m_l\leq m$, we proceed as before, i.e., we show that
\begin{equation}\label{ineq:median-lower-bound}
 \Prob{L(G)\leq m_l} \leq {1\over 2}\,.
\end{equation}
Indeed, observe that 
  since $|E'|= |E|-|E\setminus E'|$ and $L(G')\leq L(G)$,
\begin{eqnarray*}
\Prob{L(G)\leq m_l} &\leq &
    \Prob{|E|\leq \left(1-\eta\right){\asize\bsize\over k}} \\
 & & \mbox{}+ \Prob{|E\setminus E'|\geq \delta\cdot{\asize\bsize\over k}\,}\\
 & & \mbox{}+ \Prob{L(G')\leq 2(1-\delta)\sqrt{\asize\bsize\over k},
	    |E'| > (1-\eta-\delta){\asize\bsize\over k}}.
\end{eqnarray*}

We again bound the terms one by one, applying as done above
  Lemma~\ref{l:exp-del-edges}, Markov's inequality and Theorem~\ref{t:bdj}, respectively.
Indeed, for a suitable real value $\eta>0$ and 
  $\lambda=2-[2(1-\delta)/\sqrt{1-2\eta}]>0$ we get
\begin{eqnarray*}
\Prob{L(G)\leq m_l}
  \leq {1\over \eta^2 (\asize\bsize/k)} + {1\over 6}
   + B_1\exp\left(-c\lambda^{3}{\asize\bsize\over k}\right)\,.
\end{eqnarray*}
So, (\ref{ineq:median-lower-bound}) follows by taking again $\eta=\sqrt{6/C}$
and using $\asize\bsize \geq Ck$.
Proposition~\ref{p:small} is proved.
\end{proof}

\section{The lower bound in  Theorem~\ref{t:}}
In this section we establish the lower bound on the expectation
of  $L(\Sigma(K_{n,n};k))$ and the lower tail bound for its distribution.

Given $\eps$, let $\delta>0$ be such that $(1-2\delta)^{2}=1-\eps$,  
  and let $C=C(\delta)$ be as in Proposition~\ref{p:small}.
Fix $\widetilde{C}\geq \sqrt{C}$ large enough so that
\[
\exp\left(-{\delta^2\over 4(1+\delta)}\cdot \widetilde{C}\right)
  \leq \delta\,.
\]
 Let $\tilde{n}(k)=\tilde{n}=\lfloor\delta k/12\rfloor$. 
Proposition~\ref{p:small} applies for $k\geq k_0$ where $k_0$ is such that
 $\tilde{n}(k_0)\geq \widetilde{C}\sqrt{k_0}$.
It follows that
\begin{eqnarray*}
\Ex{L(\Sigma(K_{\tilde{n},\tilde{n}};k)} 
  &\geq &(1-2\delta)\cdot {2\tilde{n}\over\sqrt{k}}
  \cdot\Prob{L(G)\geq 2(1-2\delta){\tilde{n}\over\sqrt{k}}} \\
  &
  \geq & (1-2\delta)\cdot {2\tilde{n}\over\sqrt{k}}\left(1-2\exp\left(-{\delta^2\over 4(1+\delta)}\cdot {\tilde{n}\over \sqrt{k}}\right)\right)\\
&
  \geq &(1-\eps)\cdot {2\tilde{n}\over\sqrt{k}}\,.
\end{eqnarray*}
The desired lower bound on the expectation follows since by subadditivity,
  $(1/n)\cdot\Ex{L(\Sigma(K_{n,n};k)}$ is nondecreasing.

Now we establish the lower tail bound.
Let $\tilde{n}=\lceil C\sqrt{k}\,\rceil$ and 
  $q=\lfloor n/\tilde{n}\rfloor$.
Moreover, let $G$ be chosen according to $\Sigma(K_{n,n};k)$ 
  and let $G_{i}$ be the subgraph induced in $G$
  by the vertices $(i{-}1)\cdot\tilde{n}+1,\ldots,i\cdot\tilde{n}$
  in each color class, $i=1,\ldots,q$.
We observe that $L(G_1),\ldots,L(G_q)$ are independent identically
  distributed with distribution $\Sigma(K_{\tilde{n},\tilde{n}};k)$
  and $L(G)\geq L(G_1)+\cdots+L(G_q)$.
Let $\mu=\Ex{L(G_i)}$ and $t=\eps(2n/\sqrt{k})$.
Since $n\leq (q+1)\tilde{n}$,
 the lower bound on
  $\mu$ proved above yields that
\[
\Prob{L(G)\leq (1-3\eps)\cdot {2n\over \sqrt{k}}} \leq 
 \Prob{\sum_{i=1}^{q}L(G_i)\leq q\mu-t+
   (\mu-t)}.
\]
An argument similar to the one used above to derive the 
bound  $\mu\geq (1-\eps)2\tilde{n}/\sqrt{k}$ can be used to 
  obtain  
  $\mu\leq (1+\eps)2\tilde{n}/\sqrt{k}$
 from Propostion~\ref{p:small}.
Let $n$ be large enough so that $n\geq\tilde{n}(1+2\eps)/\eps$.
Thus, $q\geq (1+\eps)/\eps$
  and $t\geq \eps q\mu/(1+\eps)\geq \mu$.
Hence, 
 a standard Chernoff bound~\cite[Theorem~2.1]{jlr00} 
  implies that
\[
\Prob{L(G)\leq (1-3\eps)\cdot {2n\over \sqrt{k}}} 
 \leq
 \Prob{\sum_{i=1}^{q}L(G_i)\leq q\mu-t}
 \leq
 \exp\left(-{t^2\over 2q\mu}\right) \leq
 \exp\left(-{\eps^2\over 2(1+\eps)}\cdot {2n\over \sqrt{k}}\right).
\]

\section{The upper  bound in Theorem~\ref{t:}}
We will only discuss the tail bound since   
  $L(\Sigma(K_{n,n};k))\leq n$ always, 
  and so the claimed estimate for the expectation follows from the tail bound.

Let $\eps>0$ be fixed.
We choose a sufficiently small $\delta=\delta(\eps)>0$, much smaller
  than $\eps$.
Requirements on $\delta$ will be apparent from the subsequent proof.

Henceforth, we fix constants $1/2<\alpha<\beta<3/4$
  (any choice of $\alpha$ and $\beta$ in the specified range 
  would suffice for our purposes).
In this section, we will always assume that $k\geq k_0$ for a sufficiently
  large integer $k_0=k_0(\eps)$, and that $n$ is sufficiently
  large compared to $k$: $n\geq k^{\beta}$, say.
Note that for $n\leq k^{\beta}$ (and $k$ sufficiently large), 
  the tail bound of Theorem \ref{t:} follows
  from Proposition~\ref{p:small}.

\heading{Block partitions. }
Let us write
\[
m_{\mathrm{max}}=(1+\eps)\cdot \frac{2n} {\sqrt k}
\]
for the upper bound on the expected size of a planar matching
  as in Theorem~\ref{t:}.
We also define an auxiliary parameter
\[
\ell = k^{\alpha}\,.
\]
This is a somewhat arbitrary choice (but given by a simple formula).
The essential requirements on $\ell$ are that 
  $\ell$ be much larger than $\sqrt k$ and much smaller
  than $k^{3/4}$.
We note that $n/\ell$ is large by our assumption $n\geq k^{\beta}$.

Let $M$ be a planar matching with
  $m_{\mathrm{max}}$ edges on the sets $A$ and $B$, $|A|=|B|=n$.
We define a partition of $M$ into blocks of consecutive edges.
There will be roughly $n/\ell$ blocks, each of them containing
  at most
\[
e_{\mathrm{max}}
  = \left\lfloor{1\over\delta}\cdot{\ell\over n}\cdot
     m_{\mathrm{max}}\right\rfloor
\]
  edges of $M$.
So $e_{\mathrm{max}}$ is of order $\ell/\sqrt{k}$,
  which by our assumptions
  can be assumed to be larger than any prescribed constant.
Moreover, we require that no block is ``spread'' over more than
  $\ell$ consecutive nodes in $A$ or in $B$.

Formally, the $i$th block of the partition will be specified
  by nodes $a_i,a_i'\in A$ and $b_i,b'_i\in B$;
  $a_ib_i\in M$ is the first edge in the block and
  $a'_ib'_i\in M$ is the last edge (the block may contain only one
  edge, and so $a_ib_i=a'_ib'_i$ is possible).
The edge $a_1b_1$ is the first edge of $M$,
  and $a_{i+1}b_{i+1}$ is the edge of $M$ immediately following
  $a'_ib'_i$. Finally, given $a_ib_i$, the edge $a'_i,b'_i$
  is taken as the rightmost edge of $M$ such that
\begin{itemize}
\item the $i$th block has at most $e_{\mathrm{max}}$ edges of $M$, and
\item $a'_i-a_i\leq\ell$ and $b'_i-b_i\leq \ell$
  (here and in the sequel, with a little abuse of notation, we
  regard the nodes in $A$ and those in $B$ as natural numbers
  $1,2,\ldots,n$, although of course, the nodes in $A$ are distinct
  from those of $B$).
\end{itemize}
Let $\nblks$ denote the number of blocks obtained in this way.
It is easily seen that $\nblks=O(n/\ell)$.

\begin{figure}[h]
\begin{center}
\setlength{\unitlength}{0.00083333in}
{\renewcommand{\dashlinestretch}{30}
\begin{picture}(7661,2805)(500,250)
\put(0,0){\includegraphics{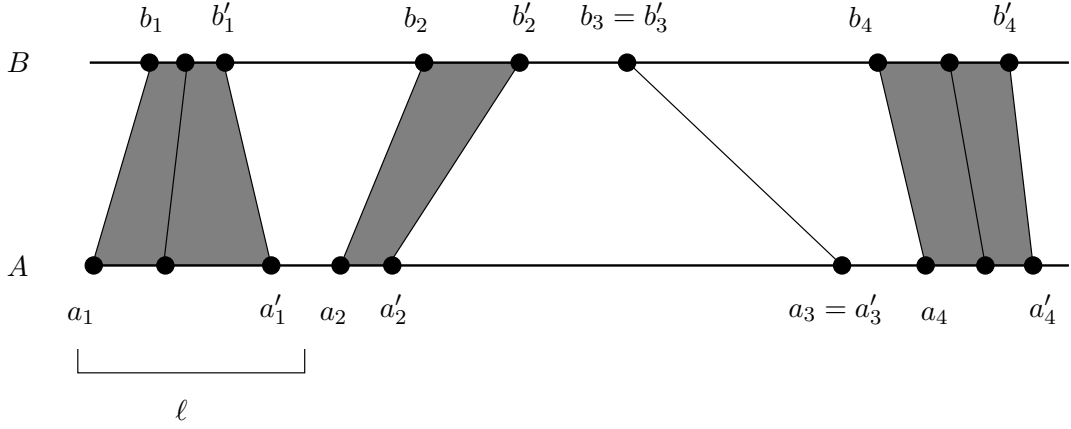}}
\put(825,975){\makebox(0,0)[lb]{$A$}}
\put(825,2250){\makebox(0,0)[lb]{$B$}}
\put(2100,2500){\makebox(0,0)[lb]{$b'_1$}}
\put(1650,2500){\makebox(0,0)[lb]{$b_1$}}
\put(3975,2500){\makebox(0,0)[lb]{$b'_2$}}
\put(3300,2500){\makebox(0,0)[lb]{$b_2$}}
\put(6975,2500){\makebox(0,0)[lb]{$b'_4$}}
\put(6075,2500){\makebox(0,0)[lb]{$b_4$}}
\put(4400,2500){\makebox(0,0)[lb]{$b_3=b'_3$}}
\put(1200,675){\makebox(0,0)[lb]{$a_1$}}
\put(2775,675){\makebox(0,0)[lb]{$a_2$}}
\put(3150,675){\makebox(0,0)[lb]{$a'_2$}}
\put(6525,675){\makebox(0,0)[lb]{$a_4$}}
\put(7200,675){\makebox(0,0)[lb]{$a'_4$}}
\put(5700,675){\makebox(0,0)[lb]{$a_3=a'_3$}}
\put(2400,675){\makebox(0,0)[lb]{$a'_1$}}
\put(1875,75){\makebox(0,0)[lb]{$\ell$}}
\end{picture}
}
\end{center}
\caption{A block partition.}\label{f:1}
\end{figure}

A block partition is schematically illustrated in Fig.~\ref{f:1}.

\heading{Counting the types. }
Let $e_i$ be the number of edges of $M$
  in the $i$th block. Let us call the $5\nblks$-tuple
  $T=(a_1,a'_1,b_1,b'_1,e_1,\ldots,a_\nblks,a'_\nblks,b_\nblks,b'_\nblks,e_\nblks)$
  the \emph{type} of the block partition of $M$, and let us write $T=T(M)$.
Let $\TT$ denote the set of all possible types of block partitions of
  planar matchings as above.

\begin{lemma}\label{l:ntypes}
We have
\[
|\TT|\leq \exp\left(C_1\frac n\ell\log \ell\right)
\]
with a suitable absolute constant $C_1$.
\end{lemma}
\begin{proof}
The number of choices for $a_1,\ldots,a_\nblks$ is at most the number
  of ways of choosing $\nblks$ elements out of $n$, i.e., $n\choose \nblks$.
Since $m_{\mathrm{max}}\leq n$, the number of choices for the $e_i$
  is no larger than the number of partitions of $n$ into
  $\nblks$ positive summands, which is $n\choose \nblks$.
Grossly overestimating, for a fixed $\nblks$ we can thus bound the number
  of types by ${n\choose \nblks}^5$.
Using the standard estimate ${n\choose \nblks}\leq
\left(en/\nblks\right)^\nblks$
  and $\nblks=O(n/\ell)$, we get $ \log |\TT|= O((n/\ell)\log\ell)$ as claimed.
\end{proof}

\heading{The probability of a matching with a given type of block partition. }
Next we show that for every fixed type $T$, the probability that our
  random graph contains a planar matching of size $m_{\mathrm{max}}$
  with that type of block partition is very small.

\begin{lemma}\label{l:pT}
Let $n$ and $k$ be as above.
For any given type $T\in\TT$, the probability $p_T$
  that the random graph $\Sigma(K_{n,n};k)$ contains a planar
  matching $M$ with $m_{\mathrm{max}}$ edges and with  $T(M)=T$
  satisfies
\[
p_T\leq \exp\left(-c\eps^2 \delta \cdot\frac n {\sqrt k}\right)
\]
with a suitable absolute constant $c>0$.
\end{lemma}
\begin{proof}
Let $G_i$ denote the subgraph of the considered random graph
$\Sigma(K_{n,n};k)$ induced by the nodes
  $a_i,a_i+1,\ldots, a'_i$ and $b_i,b_i+1,\ldots,b'_i$.
We note that the distribution of $G_i$
  is the same as that of $\Sigma(K_{\asize_i,\bsize_i};k)$,
  where $\asize_i=a'_i-a_i+1$ and $\bsize_i=b'_i-b_i+1$.

A necessary condition for the existence of a planar matching
  $M$ with $T(M)=T$ is $L(G_i)\geq e_i$ for all $i=1,2,\ldots,\nblks$.
Crucially for the proof, the events $L(G_i)\geq e_i$
  are independent for distinct $i$, and so we have
\[
p_T\leq \prod_{i=1}^\nblks \Prob{L(\Sigma(K_{\asize_i,\bsize_i};k))\geq e_i}.
\]
The plan is to apply Proposition~\ref{p:small}(i) for each $i$.
The construction of the block partition guarantees that $\asize_i,\bsize_i\leq
\ell$,
  and so  the condition $(\asize_i+\bsize_i)\sqrt{\asize_i\bsize_i}\leq 
  \delta K^{3/2}/6$ in Proposition~\ref{p:small}
  is satisfied. However,  the condition $\asize_i\bsize_i\geq Ck$ may fail.
To remedy this, we artificially enlarge the blocks; clearly, this can
  only increase the probability that a planar matching
  of size $e_i$ is present.

Let us call the $i$th block \emph{short} if it is the last block,
  i.e., $i=\nblks$, or if $e_i=e_{\mathrm{max}}$. Let $S\subseteq [\nblks]$
  denote the set of all indices of short blocks. We have
  $(|S|-1)e_{\mathrm{max}} \leq m_{\mathrm{max}}$, and since
  $e_{\mathrm{max}}\geq \frac 1\delta \cdot\frac \ell n \cdot
  m_{\mathrm{max}}-1$,
  we obtain
  $|S|\leq 2\delta n/\ell$.

The blocks that are not short are called \emph{regular},
  and we write $R=[\nblks]\setminus S$.
For a regular block $i$, we have $\max(a_{i+1}-a_i,b_{i+1}-b_i)\geq\ell$
  by the construction of the block partition.

Now we define the sizes of the artificially enlarged graphs,
which will replace the $G_i$ in the subsequent calculation.
Namely, for a short block ($i\in S$), we set
$$
\bar \asize_i = \bar \bsize_i = \ell.
$$
For a regular block ($i\in R$),
 we distinguish two cases. If $a_{i+1}-a_i\geq \ell$,
we set $\bar \asize_i =\ell$ and $ \bar \bsize_i =\max(\delta\ell,\bsize_i)$.
Otherwise, we set  $\bar \asize_i =\max(\delta\ell,\asize_i)$
and $ \bar \bsize_i =\ell$.

In the first case above, we have $\bar \asize_i\leq a_{i+1}-a_i$
and $\bar \bsize_i-\bsize_i\leq \delta\ell$, and similarly for the second
case. Therefore,
$\sum_{i\in R} \bar \asize_i\leq n+\delta \ell\cdot|R| =
(1+O(\delta))n$, with an absolute constant in the $O(\cdot)$
notation, and similarly $\sum_{i\in R} \bar \bsize_i=
(1+O(\delta))n$. For $i\in S$ we find
$\sum_{i\in S} \bar \asize_i,
\sum_{i\in S} \bar \bsize_i\leq |S|\cdot \ell \leq 2\delta n$.
Altogether
\begin{equation}\label{e:rssum}
\sum_{i=1}^\nblks \bar \asize_i \leq (1+O(\delta))n,\ \ \ \ \
\sum_{i=1}^\nblks \bar \bsize_i \leq (1+O(\delta))n.
\end{equation}

Now $\bar \asize_i$ and $\bar \bsize_i$ already satisfy the requirements
of Proposition~\ref{p:small}(i), since
we have 
$\bar \asize_i\bar \bsize_i \geq \delta\ell^2 =\delta k^{2\alpha}>Ck$
and
$(\bar\asize_i+\bar\bsize_i)\sqrt{\bar\asize_i\bar\bsize_i}\leq
2\ell^2=2k^{2\alpha} < \delta k^{3/2}/6$.
We thus have,
by Proposition~\ref{p:small},
$$
\Prob{L(\Sigma(K{\bar \asize_i,\bar \bsize_i};k))\geq e_i}
\leq 2e^{-(e_i-m_u(\bar \asize_i,\bar \bsize_i))^2/8e_i}
$$
for all $i$ such that $e_i\geq m_u(\bar \asize_i,\bar \bsize_i)$,
where $m_u(r,s)=(1+\delta) 2\sqrt{rs/k}$.
In the denominator of the exponent, we estimate
$e_i\leq e_{\mathrm{max}}$. We thus have
$$
p_T \leq \prod_{i=1}^\nblks  2e^{-\max(0,e_i-m_u(\bar \asize_i,\bar
\bsize_i))^2/8e_{\mathrm{max}}}
$$
(note that the factors for $i$ with  $e_i< m_u(\bar \asize_i,\bar \bsize_i)$
equal 1). We consider the logarithm of $p_T$, we use the
Cauchy--Schwarz inequality, and the inequality
$\max(0,x)+\max(0,y)\geq \max(0,x+y)$:
\begin{eqnarray*}
-\ln p_T &\geq &\frac 1{8e_{\mathrm{max}}}\sum_{i=1}^\nblks
\max(0,e_i-m_u(\bar \asize_i,\bar \bsize_i))^2 - q\ln 2\\
 &\geq &\frac 1{8e_{\mathrm{max}}}\cdot {1\over \nblks}
\cdot \left(\sum_{i=1}^\nblks \max(0,e_i-m_u(\bar \asize_i,\bar
\bsize_i))\right)^2 - q\ln 2\\
&\geq& \Omega(1)\cdot \frac 1{e_{\mathrm{max}}}
\cdot \frac \ell n \biggl(\sum_{i=1}^\nblks e_i -
\sum_{i=1}^\nblks m_u(\bar \asize_i,\bar \bsize_i)\biggr)^2 - q\ln 2\\
&\geq& \Omega\left(\frac {\delta \sqrt k}{n}\right) \left((1+\eps)\frac
{2n}{\sqrt k} -
\frac {2(1+\delta)}{\sqrt k} \sum_{i=1}^\nblks \sqrt{\bar \asize_i \bar
\bsize_i}\right)^2 - q\ln 2.
\end{eqnarray*}
The function $(x,y)\mapsto \sqrt{xy}$ is subadditive:
$\sqrt{xy}+\sqrt{x'y'}\leq \sqrt{(x+x')(y+y')}$.
Thus, using (\ref{e:rssum}), we have
\[
\sum_{i=1}^\nblks \sqrt{\bar \asize_i \bar \bsize_i}\leq (1+O(\delta))n\,,
\]
and so, since $q=O(n/l)$ and $l\geq\sqrt{k}$,
\[
-\ln p_T 
  \geq \Omega\left({\delta\sqrt k\over n}\right)\left((1+\eps){2n\over\sqrt k}
    - (1+O(\delta)){2n\over\sqrt k}\right)^2 - q\ln 2
  =\Omega\left(\eps^2\delta\cdot{n\over\sqrt k}\right)\,.
\] 
Lemma~\ref{l:pT} is proved.
\end{proof}

\heading{Proof of Theorem \ref{t:}. }
We have
\[
\Prob{L(\Sigma(K_{n,n};k))\geq m_{\mathrm{max}}} \leq
\sum_{T\in \TT} p_T \leq |\TT|\cdot \max_T p_T\,.
\]
The sought after estimate
\[
\Prob{L(\Sigma(K_{n,n};k))\geq
m_{\mathrm{max}}}\leq \exp\left(-\Omega(\eps^2\delta n/\sqrt k\,)\right)\,,
\]
follows from Lemmas~\ref{l:ntypes} and~\ref{l:pT}.
\proofend

\section{Extensions}\label{s:bin}
Similarly one can prove results  for the Erd\H{o}s model 
  analogous to those obtained in previous sections
  (essentially, $k$ is now replaced by $1/p$):
\begin{theorem}\label{t:bin}
For every $\eps>0$
  there exist constants $p_0\in (0,1)$ and $C$ such that for all
  $p<p_0$ and all $n$ with $n\sqrt{p}> C$ we have
\[
  (1-\eps)\cdot 2n\cdot\sqrt{p}\leq
  \Ex{L(G(K_{n,n};p))} \leq 
  (1+\eps)\cdot 2n\cdot\sqrt{p}\,.
\]
Moreover, there is an exponentially small tail bound; namely, for every
  $\eps>0$ there exists $c>0$ such that for $p$ and $n$ as above,
$$
\Prob{\left| L(G(K_{n,n};p)) - 2n\sqrt{p}\right|\geq \eps 2n\sqrt{p}}
        \leq e^{-cn\sqrt{p}}\,.
$$
\end{theorem}
Subadditivity arguments yield that $\Ex{L(G(K_{n,n};p))}/n$ converges 
  to a constant $\Delta_p$ as $n\to\infty$.
The previous theorem thus implies that $\Delta_p/\sqrt{p}\to 2$ as
  $p\to 0$.

Also, similar results hold for the $G(K_{r,s};p)$ model as those derived 
  for $\Sigma(K_{\asize,\bsize};k)$.
Specifically,
\begin{prop}\label{p:small-bin}
For every $\delta>0$, there exists a (large) positive constant $C$
such that:

\begin{enumerate}
\item[\rm (i)]
If $\asize\bsize\geq C/p$ and
$(\asize+\bsize)\sqrt{\asize\bsize}\leq \delta/6 p^{3/2}$,
then  with $m_u=m_u(\asize,\bsize)=2(1+\delta)\sqrt{\asize\bsize p}$, we have
\[
 \Prob{L(G(K_{\asize,\bsize};p))\geq m_u+t}\leq 2e^{-t^2/8(m_u+t)}
\]
for all $t\geq 0$.
\item[\rm (ii)]
If $\asize\bsize\geq C/p$ and
and $\asize+\bsize\leq \delta/6p$, then with
$m_u$ as above and $m_l=m_l(\asize,\bsize)=2(1-\delta)\sqrt{\asize\bsize p}$,
we have
\[
 \Prob{L(G(K_{\asize,\bsize};p))\leq m_l-t}\leq 2e^{-t^2/8m_u}
\]
for all $t\geq 0$.
\end{enumerate}
\end{prop}

In~\cite{johansson00}, Johansson implicitly considers a model somewhat
  related to the $G(K_{n,n};p)$ model.
Specifically, a distribution $G^{*}(K_{n,n};p)$ over weighted
  instances of $K_{n,n}$.
The weight of each edge is a geometrically distributed random
  variable taking the value $k\in\N$ with probability 
  $(1-p)^{k}p$, and the edge weights are mutually independent.
Denoting the maximum weight planar matching  
  of an instance drawn according to $G^{*}(K_{n,n};p)$ 
  by $L(G^{*}(K_{n,n};p))$, Johansson's result~\cite[Theorem 1.1]{johansson00}
  says that for all $pi\in (0,1)$,
\begin{eqnarray*}
\lim_{n\to\infty} {1\over n}\cdot\Ex{L(G^{*}(K_{n,n};p))} =
  {(1+\sqrt{1-p})^{2}\over p}\,.
\end{eqnarray*}
Note that an instance $G$ of $G(K_{n,n};p)$ can be obtained 
  from one drawn according to $G^{*}(K_{n,n};p)$ by including in $G$
  only those edges of $K_{n,n}$ with nonzero weight.
Hence, 
\begin{eqnarray*}
\Ex{L(G(K_{n,n};p))} & \leq & \Ex{L(G^{*}(K_{n,n};p))}\,.
\end{eqnarray*}
It follows that $\Delta_{p} \leq {(1+\sqrt{1-p})^{2}/p}$,
  for all $p\in (0,1)$.
We shall see below that known results imply a much stronger
  bound on $\Delta_p$ for not too large values of $p$.

Gravner, Tracy and Widom~\cite{gtw01} consider processes 
  associated to random $(0,1)$--matrices 
  where each entry takes the value $1$ with probability
  $p$, independent of the values of other matrix entries.
 In particular they study a process called
  \emph{oriented digital boiling} (ODB) and analyze the behavior
  of a so called \emph{height function} which equals, in 
  distribution, the longest sequence $(i_l,j_l)$ of 
  positions in a random $(0,1)$--matrix of size $n\times n$
  which have entry $1$ such that the $i_l$'s are increasing
  and the $j_l$'s are nondecreasing.
In contrast, $L(G(K_{n,n};p))$ equals in distribution the 
  longest such sequence with both $i_l$'s and $j_l$'s increasing.
This latter model is referred to 
  as \emph{strict oriented digital boiling} in~\cite{gtw01},
but no results are claimed for it.
Clearly, an ODB process dominates that of a strict ODB process. 
Hence, \cite[\S 3, (1)]{gtw01} implies that for any $p<1/2$,
\[
\Delta_p \leq \kappa_p:=
  \lim_{n\to\infty} {1\over n}\cdot\Ex{L(G(K_{n,n};p))} = 2\sqrt{p(1-p)}\,,
\]
which in turn implies that $\limsup_{p\to 0} 
  \Delta_{p}/\sqrt{p}\leq 2$.
Nevertheless, our derivation of this latter limit value is elementary
  in comparison with the highly technical nature of~\cite{gtw01}.

\section*{Acknowledgments}
We thank Ricardo Baeza for calling to our attention reference~\cite{sm83}.

\bibliographystyle{plain}
\bibliography{match}

\end{document}